\documentclass[12pt]{article}

\usepackage{a4}   
\usepackage{amssymb,amsmath,mathrsfs,textcomp}   
\usepackage[amsmath,thmmarks]{ntheorem}   
\usepackage{yhmath}   
\usepackage{indentfirst}   
\usepackage{tocloft}   
\usepackage[hyperfootnotes=false,colorlinks=true,linkcolor=black,citecolor=black]{hyperref}   

\usepackage{graphicx}
\usepackage{epstopdf}

\theoremstyle{nonumberplain}
\theorembodyfont{\rm}
\theoremseparator{.}

\theoremstyle{plain}
\theorembodyfont{\rm}
\theoremseparator{.}

\theoremstyle{plain}
\theorembodyfont{\slshape}
\theoremseparator{.}

\theoremstyle{plain}
\theorembodyfont{\slshape}
\theoremseparator{$^\prime$.}

\theoremstyle{nonumberplain}
\theorembodyfont{\slshape}
\theoremseparator{.}

\theoremstyle{nonumberplain}
\theoremheaderfont{\bf}
\theorembodyfont{\rm} 
\theoremseparator{.}
\theoremsymbol{$\Box$}
\newtheorem{proof}{Proof}

\makeatletter
\renewcommand*{\@seccntformat}[1]{\csname the#1\endcsname.\quad}
\makeatother

\newcommand{\coloneqq}{\mathrel{\mathop:}=}
\newcommand{\eqqcolon}{=\mathrel{\mathop:}}

\renewcommand{\Re}{\text{Re}\,}
\renewcommand{\Im}{\text{Im}\,}

\makeatletter
\def\moverlay{\mathpalette\mov@rlay}
\def\mov@rlay#1#2{\leavevmode\vtop{%
   \baselineskip\z@skip \lineskiplimit-\maxdimen
   \ialign{\hfil$\m@th#1##$\hfil\cr#2\crcr}}}
\newcommand{\charfusion}[3][\mathord]{
    #1{\ifx#1\mathop\vphantom{#2}\fi
        \mathpalette\mov@rlay{#2\cr#3}
      }
    \ifx#1\mathop\expandafter\displaylimits\fi}
\makeatother

\begin{document}

%
%
\title{On the existence of Kobayashi and Bergman metrics for Model domains}
\author{Nikolay Shcherbina}
\maketitle
\nopagebreak

%
%
\small\noindent{\bf Abstract.}
We prove that for a pseudoconvex domain of the form $\mathfrak{A} = \{(z, w) \in \mathbb C^2 : v > F(z, u)\}$, where $w = u + iv$ and F is a continuous function on ${\mathbb C}_z \times {\mathbb R}_u$, the following conditions are equivalent:
\begin{enumerate}
  \item[(1)] The domain $\mathfrak{A}$ is Kobayashi hyperbolic.
  
  \item[(2)] The domain $\mathfrak{A}$ is Brody hyperbolic.  

  \item[(3)] The domain $\mathfrak{A}$ possesses a Bergman metric.

  \item[(4)] The domain $\mathfrak{A}$ possesses a bounded smooth strictly plurisubharmonic function, i.e. the core $\mathfrak{c}(\mathfrak{A})$ of $\mathfrak{A}$ is empty.

  \item[(5)]  The graph $\Gamma(F)$ of $F$ can not be represented as a foliation by holomorphic curves of a very special form, namely, as a foliation by translations of the graph $\Gamma({\mathcal H})$ of just one entire function ${\mathcal H} : {\mathbb C}_z \to {\mathbb C}_w$. \normalsize
\end{enumerate}

%
%
\renewcommand{\thefootnote}{}

%
%
\setcounter{tocdepth}{1}
\tableofcontents

\section{Introduction}

\noindent
The purpose of this paper is to study the Kobayashi and Bergman metrics for pseudoconvex domains of the form $\mathfrak{A} = \{(z, w) \in {\mathbb C^n_z \times {\mathbb C}_w} : v > F(z, u)\}$, where $w = u + iv$ and F is a continuous function on $\mathbb C^n_z \times {\mathbb R}_u$. This type of domains we call in what follows for {\it Model domains}. They appear naturally (usually with much more special choice of the function F) as the limit domains in the scaling method (see, for example, \cite{Pinchuk91}), in the representations of domains of finite type, in the biholomorphic classification problem etc. Note that in the last problem the existence of Kobayashi metric on $\mathfrak{A}$ is of special interest, since in the case of its existence the group ${\text{Aut}}(\mathfrak{A})$ of holomorphic automorphisms of $\mathfrak{A}$ is in fact a finite-dimensional real Lie group (see, for example, [Ka, Folgerung 2.6, p. 55] or [Ko1, Theorem 2.1, p. 68]).

Despite the fact that many properties for special classes of such domains were intensively studied (see, just to cite a few, \cite{KohnNirenberg73}, \cite{Fornaess77}, \cite{BedfordFornaess78}, \cite{BedfordPinchuk94}, \cite{CoupetPinchuk01}, \cite{ChenKamimotoOhsawa04}, \cite{AhnGaussierKim16}), to the best of our knowledge no results were known in the general case. Surprisingly, in the case of complex dimension two there is a complete characterization of Model domains possessing Kobayasi and Bergman metrics in the general setting.

The main result of this paper asserts that for a Model domain $\mathfrak{A}$ in $\mathbb C^2$ the property of being Kobayashi hyperbolic and the property to possess a Bergman metric occur simultaneously. Moreover, these two properties hold true for all such domains, except for the very specal case when the boundary $\partial{\mathfrak{A}}$ of $\mathfrak{A}$ is foliated by translations of just one entire function $w = {\mathcal H}(z)$. More precisely, the following statement holds true.

\vspace{0,3cm}
\noindent \textbf{Main Theorem.} {\it Let $\mathfrak{A}$ be a pseudoconvex domain of the form 
\[ \mathfrak{A} = \{(z, w) \in \mathbb C^2 : v > F(z, u)\}, \]
where $w = u + iv$ and F is a continuous function on ${\mathbb C}_z \times {\mathbb R}_u$. Then the following conditions are equivalent:}
\begin{enumerate}
  \item[(1)] The domain $\mathfrak{A}$ is Kobayashi hyperbolic.
  
  \item[(2)] The domain $\mathfrak{A}$ is Brody hyperbolic.  

  \item[(3)] The domain $\mathfrak{A}$ possesses a Bergman metric.

  \item[(4)] The domain $\mathfrak{A}$ possesses a bounded smooth strictly plurisubharmonic function, i.e. the core $\mathfrak{c}(\mathfrak{A})$ of $\mathfrak{A}$ is empty.

  \item[(5)]  {\it The boundary $\partial{\mathfrak{A}}$ of $\mathfrak{A}$ can not be represented as follows:
\[ \partial{\mathfrak{A}} = \bigcup_{t \in {\mathbb R}} \Gamma({\mathcal H} + (t +i\Theta(t))), \]
where ${\mathcal H} : {\mathbb C}_z \to {\mathbb C}_w$  is an entire function, $\Theta : {\mathbb R} \to {\mathbb R}$ is a continuous function and $\Gamma({\mathcal H} + (t +i\Theta(t)))$ denote the graph of the function $w = {\mathcal H}(z) + (t +i\Theta(t))$.}
\end{enumerate}

\vspace{0,3cm}
Condition (4) of the theorem above plays a crucial role in its proof and can also be expressed in terms of the core $\mathfrak{c}(\mathfrak{A})$ of $\mathfrak{A}$ (for the definition of the core see the next section of this paper). The notion of the core of a domain (or, more generally, the core of a manifold) was introduced and systematically studied few years ago by Harz-Shcherbina-Tomassini in [HST1-3]. Some results of [HST2] on the structure of the core were then generalised to the case of higher dimensions by Poletsky-Shcherbina in \cite{PoletskyShcherbina19} and Slodkowski in \cite{Slodkowski19}. Recently Poletsky in \cite{Poletsky19} used manifolds with empty core to develop further the theory of pluricomplex Green functions. A notion similar to the notion of a manifold with empty core, but with no regularity assumptions on plurisubharmonic functions in question, was considered by Boucksom-Diverio in \cite{BoucksomDiverio18} in connection to Lang's conjecture. They call such manifolds for manifolds of {\it bounded type}. We do not know if for domains in $\mathbb C^n$ with not "too irregular" boundary the properties of having an empty core and being of bounded type coincide or not.

Note also that even though the notion of the core of a domain depends on the smoothness of the considered class of plurisubharmonic functions, the classical Richberg's theorem (see, for example, Theorem I.5.21 in \cite{Demailly12} or \cite{Richberg68}) implies that the condition of having an empty core does not depend on the smoothness class as soon as the functions are at least continuous.

The proof of our theorem is essentially based on the results and technique developed in \cite{Shcherbina93}. For the convenience of the reader we will remind in the next section the main result and the necessary preparatory statements of that paper.

As a final comment, we want to point out that the results of our Main Theorem are not true in general for Model domains in $\mathbb C^n$ with $n \geq 3$ as it is shown in the forthcoming paper \cite{GaussierShcherbina19}.

\section{Preliminaries} \label{sec_preliminaries}

\noindent
In this section we collect some definitions and results which will be used in the proof of the Main Theorem.

We start with the definition of the core of a manifold.

\vspace{0,3cm}
\noindent \textbf{Definition 1.} Let $\mathcal{M}$ be a complex manifold. Then the set
\begin{equation*}\begin{split}
\mathfrak{c}(\mathcal{M}) \coloneqq \big\{&p \in \mathcal{M} : \text{every smooth plurisubharmonic function on $\mathcal{M}$ that is} \\& \text{bounded from above fails to be strictly plurisubharmonic in } p \big\}
\end{split}\end{equation*}
is called the \emph{core} of $\mathcal{M}$.

\vspace{3mm}
The next statement, which follows easily from the definition of the core, shows that two properties mentioned in the Condition (4) of the Main Theorem are equivalent.

\vspace{3mm}
\noindent \textbf{Lemma 1.} {\it For a complex manifold $M$ the property that the core $\mathfrak{c}(M)$ of $M$ is empty holds true if and only if there is a bounded smooth strictly plurisubharmonic function on $M$.}

\begin{proof} 
Assume that the core $\mathfrak{c}(M)$ of $M$ is empty. Then for every $p \in M$, there exists a smooth bounded above plurisubharmonic function $\psi_p$ on $M$ that is strictly plurisubharmonic on an open neighbourhood $V_p \subset M$ of $p$. Let $\{p_j\}_{j=1}^\infty$ be a sequence of points $p_j \in M$ such that $\bigcup_{j=1}^\infty V_{p_j} = M$. Then, for suitably chosen sequence of positive numbers $\{\varepsilon_j\}_{j=1}^\infty$ the function $\tilde{\psi} = \sum_{j=1}^\infty \varepsilon_j \psi_{p_j}$ is a smooth bounded above function which is strictly plurisubharmonic on the whole of $M$. Hence $e^{\tilde{\psi}}$ will be a function as desired. The opposite implication is obvious.
\end{proof}

Now we formulate a shortened version of the Main Theorem in \cite{Shcherbina93} which is enough for our purposes here.

\vspace{3mm}
\noindent \textbf{Theorem 1.} {\it Let $G$ be a bounded strictly convex domain in ${\mathbb C}_z \times {\mathbb R}_u$ and $\varphi : \partial{G} \to {\mathbb R}_v$ be a continuous function. Then the following properties hold:}
\begin{enumerate}
  \item[(1)] {\it  There is a continuous function $\Phi : \bar{G} \rightarrow {\mathbb R}_v$ such that $\Phi|_{\partial{G}} = \varphi$ and $\widehat{\Gamma}(\varphi) = {\Gamma}(\Phi)$.}

  \item[(2)] {\it  The set ${\Gamma}(\Phi) \setminus \Gamma(\varphi)$ is the union of a disjoint family of complex disks $\{D_\alpha\}$.}

  \item[(3)] {\it   For each $\alpha$, there is a simply-connected domain $\Omega_\alpha \subset \mathbb C_z$ and a holomorphic function $w = {f_\alpha}(z)$, defined on the domain $\Omega_\alpha$, such that the disk $D_\alpha$ is the graph of the function $f_\alpha$.}
\end{enumerate}

\vspace{0,1cm}
\noindent
{\it Here $\Gamma(\varphi)$ means the graph of the function $\varphi$ and $\widehat{\Gamma}(\varphi)$ the polynomial hull of this graph.}

\vspace{4mm}
We also state below Lemma 3.2, Lemma 3.3 and slightly strengthened version of Lemma 3.5 from the paper \cite{Shcherbina93} enumerating them as Lemma 2, Lemma 3 and Lemma 4, respectively.

\vspace{3mm}
\noindent \textbf{Lemma 2.} {\it Let $\mathcal U$ be a smooth connected surface which is properly embedded into some convex domain $G \subset {\mathbb C}_z \times {\mathbb R}_u$. Suppose that near each point of this surface it can be defined locally by the equation $u = u(z)$. Then the surface $\mathcal U$ can be represented globally as a graph of some function $u = U(z)$, defined on some domain $\Omega \subset {\mathbb C}_z$.}

\vspace{3mm}
\noindent \textbf{Lemma 3.} {\it Let $\Omega$ be a domain in ${\mathbb C}_z$, and $U$ be a harmonic function on this domain. Suppose that the graph of the function $U$ is properly embedded into some convex domain $ G \subset {\mathbb C}_z \times {\mathbb R}_u$. If $U$ has a single-valued harmonic conjugate function $V$, then $\Omega$ is simply-connected.}

\vspace{3mm}
\noindent \textbf{Lemma 4.} {\it Let $\Omega \subset \Delta_R(0) = \{ z : |z | < R \}$ be a simply-connected domain which contains the origin. Let $U$ be a harmonic function on $\Omega$ such that its graph is properly embedded into the cylinder $\Delta_R(0) \times (-R, R) \subset {\mathbb C}_z \times {\mathbb R}_u$ and $U(0) = 0$. Suppose that $U$ has  a harmonic conjugate function $V$ such that $\inf_{z \in \Omega}V(z) \geq M$ and $V(0) \leq L$. Then there exists $R^\ast > 0$, which depends only on $R$, $M$ and $L$, such that the disk $\Delta_{R^\ast}(0)$ is contained in $\Omega$.}
 
The proof of Lemma 4 here goes exactly the same way as the proof of Lemma 3.5 in \cite{Shcherbina93}. We only need to replace there the rectangle $\Pi_{R, M} = \{w \in {\mathbb C}_w : -R < u < R, -M < v < M\}$ with the half-strip $\widetilde{\Pi}_{R, M} = \{w \in {\mathbb C}_w : -R < u < R, -M < v \}$, the set $E = \bar{\Pi}_{R, M} \cap \{v = \pm M \}$ with the set $\widetilde{E} = \bar{\Pi}_{R, M} \cap \{v = - M \}$, and the harmonic measure $\omega (0, E, \Pi_{R, M})$ with the harmonic measure $\omega (iL, \widetilde{E}, \widetilde{\Pi}_{R, M})$.

\vspace{3mm}
Now we recall some results on Kobayashi and Bergman metrics which will be used in the proof of the Main Theorem. For background on these metrics we refer to \cite{Kobayashi70}, \cite{Kobayashi98} and \cite{GreeneKimKrantz11}.

\vspace{3mm}
The next statement gives a link between Kobayashi hyperbolicity and the existence of bounded strictly plurisubharmonic functions. It is just a slight reformulation of Theorem 3 on p. 362 in \cite{Sibony81}.

\vspace{3mm}
\noindent \textbf{Theorem 2.} {\it Let $M$ be a complex manifold which has a bounded continuous strictly plurisubharmonic function. Then $M$ is Kobayashi hyperbolic.}

\vspace{3mm}
A similar criterion holds true for the existence of Bergman metric. It is a consequence of Theorem 1, p. 2998, and Observation 2, p. 3002, in \cite{ChenZhang02}.

\vspace{3mm}
\noindent \textbf{Theorem 3.} {\it Let $M$ be a Stein manifold which has a bounded continuous strictly plurisubharmonic function. Then $M$ possesses a Bergman metric.} 

\vspace{3mm}
We also remind here the definition of Brody hyperbolicity.

\vspace{0,3cm}
\noindent \textbf{Definition 2.} A complex manifold is said to be Brody hyperbolic if it admits no nonconstant holomorphic maps from $\mathbb C$.

\vspace{3mm}
The last preparatory statements which are needed for our proof provide a link between polynomial convexity and globally defined plurisubharmonic functions.

\vspace{3mm}
Recall first the definition of polynomial convexity.

\vspace{0,3cm}
\noindent \textbf{Definition 3.} For a compact set $K$ in ${\mathbb C}^n$, the {\it polynomial hull} $\widehat{K}$ of $K$ is defined as 
\[  \widehat{K} = \{p \in {\mathbb C}^n : |{\mathcal P}(p)| \leq \sup_{q \in K}|{\mathcal P}(q)| \,\, {\rm for \,\, all \,\, holomorphic \,\, polynomials} \,\, \mathcal P \,\, {\rm in} \,\, {\mathbb C}^n\}. \]
The set $K$ is called {\it polynomially convex} if $\widehat{K} = K$.

\vspace{3mm}
The following simple lemma is classical and follows, for example, from Theorem 4.3.4 in \cite{Hormander90}.

\vspace{3mm}
\noindent \textbf{Lemma 5.} {\it  Let $K$ be a compact set in ${\mathbb C}^n$. Then $p \in {\mathbb C}^n \setminus \widehat{K}$ if and only if there is a function $\phi$, plurisubharmonic in ${\mathbb C}^n$, such that}
\begin{equation} \label{polconv}
\sup_{q \in K}{\phi(q)} < \phi(p).
\end{equation}

The next statement is more involved. It is a slightly shortened version of Theorem 1.3.8 on p. 24 in \cite{Stout07}.

\vspace{3mm}
\noindent \textbf{Theorem 4.} {\it A compact set $K \subset {\mathbb C}^n$ is polynomially convex if and only if there is a non-negative smooth plurisubharmonic function $\phi$ defined on the whole of ${\mathbb C}^n$ such that $K = \{p \in {\mathbb C}^n : \phi(p) = 0 \}$ and $\phi$ is strictly plurisubharmonic on ${\mathbb C}^n \setminus K$.}

\section{Construction of the domain $\widetilde{\mathfrak{A}}$}

\noindent
For an arbitrary $R > 0$ we consider the ball $B_R(0) = \{ (z, u) \in {\mathbb C}_z \times {\mathbb R_u} : |z|^2 + u^2 < R^2 \}$ of radius $R$ centered at the origin and then denote by $\varphi_R$ the restriction $F|_{\partial{B_R(0)}}$ of the function $F$ to the boundary ${\partial{B_R(0)}}$ of this ball. If we apply now Part (1) of the Theorem 1 above to the ball $B_R(0)$ on the place of $G$ and the function $\varphi_R$ on the place of $\varphi$, we will get a continuous extension $\Phi_R$ of the function $\varphi_R$ to $\overline{{B}_R(0)}$ whose graph is Levi flat over $B_R(0)$.

\vspace{3mm}
The following properties of the function $\Phi_R$ are easy to prove.

\vspace{3mm}
\noindent \textbf{Lemma 6.} {\it For every $R > 0$ and each $(z, u) \in B_R(0)$ one has $F(z, u) \leq \Phi_R(z, u)$.}
 
\begin{proof} 
We argue by contradiction and suppose that $F(z_0, u_0) > \Phi_R(z_0, u_0)$ for some $(z_0, u_0) \in B_R(0)$. Then
\[ t^\ast  \coloneqq \inf\{ t \in (0, \infty) : F(z, u) \leq \Phi_R(z, u) + t {\rm \,\, for \,\, all \,\,} (z, u) \in B_R(0) \} > 0, \]
and hence, by continuity of the functions $F$ and $\Phi_R$, there is a point $(z^\ast, u^\ast) \in B_R(0)$ such that $F(z^\ast, u^\ast) = \Phi_R(z^\ast, u^\ast) + t^\ast \eqqcolon v^\ast$. Since, by Part (2) of Theorem 1, the graph of the function $\Phi_R$ is foliated by holomorphic disks $\{D_\alpha\}$, there is a disk $D_{\alpha^\ast}$ such that $p^\ast \coloneqq (z^\ast, u^\ast + iv^\ast) \in D_{\alpha^\ast}$. For each $t \in [0, 1]$ we denote by $D^t_{\alpha^\ast}$ a translation of the disk $D_{\alpha^\ast}$ on $t$ in the {\it v}-direction. Then, by construction, the family $\{D^t_{\alpha^\ast}\}_{t \in [0, 1]}$ has the properties that $D^t_{\alpha^\ast} \subset \mathfrak{A} = \{(z, w) \in \mathbb C^2 : v > F(z, u)\}$ for all $t > 0$ and that $p^\ast \in D^0_{\alpha^\ast} \cap \partial{\mathfrak{A}}$. Thus, by Kontinuit{\"a}tssatz, the domain $\mathfrak{A}$ is not pseudoconvex, which contadicts to our assumptions on $\mathfrak{A}$.
\end{proof}

\noindent 
\textbf{Lemma 7.} {\it For all \, $0 < R_1 < R_2$ and each $(z, u) \in B_{R_1}(0)$ one has \, $\Phi_{R_1}(z, u) \leq \Phi_{R_2}(z, u)$.}

\begin{proof} 
First, we define for each $R > 0$ a function $F_R : {\mathbb C}_z \times {\mathbb R_u} \to \mathbb R_v$ as
\begin{equation} \label{F_R}
F_R(z, u) \coloneqq \left\{ \begin{array}{c@{\,,\quad}l} F(z, u)  &  {\rm{for}} \,\,  (z, u) \in ({\mathbb C}_z \times {\mathbb R_u}) \setminus \overline{B_R(0)}, \\ \Phi_R(z, u) &   {\rm{for}} \,\,  (z, u) \in B_R(0).   \end{array} \right.
\end{equation}
It follows from continuity of the functions $F$ and $\Phi_R$  that the function $F_R$ is continuous. Moreover, from pseudoconvexity of the domain $\mathfrak{A} = \{(z, w) : v > F(z, u)\}$ and Levi flatness of the graph $\Gamma(\Phi_R) = \{(z, w) \in B_R(0) \times {\mathbb R_v} : v = \Phi_R(z, u)\}$ of $\Phi_R$ we conclude that the domain
\begin{equation} \label{Omega_R} 
\mathfrak{A}_R \coloneqq \{(z, w) \in \mathbb C^2 : v > F_R(z, u)\}
\end{equation}
is also pseudoconvex. The statement of the lemma is now a direct consequence of Lemma 6, if we replace there the function $F$ which defines $\mathfrak{A}$ with the new function $F_{R_1}$, the domain $\mathfrak{A}$ with the domain $\mathfrak{A}_{R_1}$ and the radius $R$ with $R_2$. 
\end{proof}

Next, for each $R > 0$ we define the domain ${\widetilde{\mathfrak{A}}}_R$ as
\[ {\widetilde{\mathfrak{A}}}_R \coloneqq \{(z, w) \in B_R(0) \times \mathbb R_v : F(z, u) < v < \Phi_R(z, u)\}. \]

Observe that one of the direct consequences of Lemma 6 and Lemma 7 is the following property of the domains ${\widetilde{\mathfrak{A}}}_R$.

\vspace{3mm}
\noindent \textbf{Lemma 8.} {\it For all \, $0 < R_1 < R_2$ one has ${\widetilde{\mathfrak{A}}}_{R_1} \subset \, {\widetilde{\mathfrak{A}}}_{R_2} \subset \, \mathfrak{A}$.}

\vspace{3mm}
Now we can finally define the domain $\widetilde{\mathfrak{A}}$ as
\[ \widetilde{\mathfrak{A}} \coloneqq \bigcup_{n \in \mathbb{N}}{\widetilde{\mathfrak{A}}}_n.  \]

$\widetilde{\mathfrak{A}}$ is obviously a subdomain (which can also be empty) of ${\mathfrak{A}}$. In the next two sections we consider two possible cases for $\widetilde{\mathfrak{A}}$.

\section{The case when $\widetilde{\mathfrak{A}} = \mathfrak{A}$}

\noindent
Let us consider an arbitrary number $n \in \mathbb N$ which we fix till the end of the proof of Lemma 9. Let
\[ C_n = \max\{F_n(z, u), (z, u) \in \overline{{B}_{n + 1}(0)}\}, \]
where $F_n$ is the function $F_R$ defined by formula (2) above with $R = n$. Consider the set
\[ K_n = \{(z, w) \in \overline{{B}_{n + 1}(0)} \times \mathbb R_v : F_n(z, u) \leq v \leq C_n + 1 \}. \]

\noindent 
\textbf{Lemma 9.} {\it The set $K_n$ is polynomially convex.}

\begin{proof} 
We consider two sets
\[ A_1 = \{(z, w) \in {\mathbb C}^2 : |z|^2 + |u|^2 >(n + 1)^2 \,\, {\rm or} \,\, v > C_n + 1 \} \]
and
\[ A_2 = \{(z, w) \in \overline{{B}_{n + 1}(0)} \times \mathbb R_v : v < F_n(z, u)\}. \]
Since $A_1 \cup A_2 = \mathbb C^2 \setminus K_n$, we conclude that for proving the polynomial convexity of $K_n$ it is enough to show that $A_1 \cap {\widehat K}_n = \emptyset$ and $A_2 \cap {\widehat K}_n = \emptyset$. 

Let $(z_0, w_0)$ be a point of the set $A_1$. Then the inequality (1) will be satisfied for the point $p = (z_0, w_0)$, the set $K_n$ and the function
\[ \phi(z, w) = \max\{|z|^2 + |u|^2 - (n + 1)^2, v - (C_n + 1) \}. \]
Hence, by Lemma 5, $(z_0, w_0)$ can not be a point of the set ${\widehat K}_n$, i.e. $A_1 \cap {\widehat K}_n = \emptyset$.

Let now $(z_0, w_0)$ be a point of the set $A_2$ and assume, to get a contradiction, that $(z_0, w_0) \in \widehat K_n$. Then we consider a 1-parameter family of domains $\{\mathfrak{A}^t_n\}_{t \in [0, \infty)}$ defined as
\[ \mathfrak{A}^t_n \coloneqq \{(z, w) \in \mathbb C^2 : v > F_n(z, u) - t \}, \]
and notice that the domain $\mathfrak{A}^0_n$ of this family coincides with the defined by (3) domain $\mathfrak{A}_R$ for $R = n$. Then, since $(z_0, w_0) \in A_2 \cap \widehat K_n$, we conclude that
\[ t^\ast  \coloneqq \inf\{ t \in [0, \infty) : \widehat K_n \subset \overline{\mathfrak{A}^t_n}\} > 0. \]
Observe now that, in view of pseudoconvexity of the domain $\mathfrak{A}^0_n$, there is a smooth strictly plurisubharmonic function $\tau$ defined on $\mathfrak{A}^0_n$ such that $\tau(z, w) \to + \infty$ as $(z, w) \to \partial{\mathfrak{A}^0_n}$. Then for each $t \in [0, + \infty)$ the function $\tau^t(z, w) \coloneqq \tau(z, w + it)$ is defined on the domain $\mathfrak{A}^t_n$ and has there the same properties as the function $\tau$ on $\mathfrak{A}^0_n$. It follows from the definition of $t^\ast$ and the properties of the function $\tau$ that 
\[ \sup_{p \in {\widehat K}_n}\tau^{t^\ast + \varepsilon}(p) \to + \infty  \]
when $\varepsilon \searrow 0$. Since $t^\ast > 0$, and since $K_n \subset \overline{\mathfrak{A}^0_n}$, it also follows from the definition of the domain $\mathfrak{A}^t_n$ and the function $\tau^t$ that
\[ \sup_{p \in K_n}\tau^{t^\ast + \varepsilon}(p) \]
stays uniformly bounded as $\varepsilon \searrow 0$. Then, in view of Sard's theorem, we can choose $\varepsilon^\ast > 0$ generic so that for the value
\[ C^\ast \coloneqq \sup_{p \in {\widehat K}_n}\tau^{t^\ast + \varepsilon^\ast}(p) \]
the level set 
\[ \mathfrak{M}^\ast \coloneqq \{ p \in \mathfrak{A}^{t^\ast + \varepsilon^\ast}_n  : \tau^{t^\ast + \varepsilon^\ast}(p) = C^\ast \} \]
of the function $\tau^{t^\ast + \varepsilon^\ast}$ will be smooth and, moreover, if we choose $\varepsilon^\ast$ small enough, then for every point $p^\ast \in \mathfrak{M}^\ast \cap {\widehat K}_n$ we will also have that 
\[ \sup_{p \in K_n}\tau^{t^\ast + \varepsilon^\ast}(p)  < \tau^{t^\ast + \varepsilon^\ast}(p^\ast). \]
Observe that the last inequality implies that
\begin{equation} \label{p_ast} 
p^\ast \in {\widehat K}_n \setminus K_n.
\end{equation}
On the other hand, since the domain 
\[ \mathfrak{A}^\ast_n \coloneqq  \{ p \in \mathfrak{A}^{t^\ast + \varepsilon^\ast}_n   : \tau^{t^\ast + \varepsilon^\ast}(p) < C^\ast \} \]
is strictly pseudoconvex and has a smooth boundary, and since ${\widehat K}_n \subset \overline{\mathfrak{A}^\ast_n}$ and $p^\ast \in \partial{\mathfrak{A}^\ast_n} \cap {\widehat K}_n$, we conclude that $p^\ast$ is the local peak point for the algebra $\mathfrak{P}$ of uniform limits on ${\widehat K}_n$ of holomorphic polynomials. Then, by Rossi's "local maximum modulus principle" (see, for example, [G, Theorem 8.2, p. 92] or \cite{Rossi60}), we conclude that $p^\ast$ is a global peak point for the algebra $\mathfrak{P}$. This contradicts (4), and proves Lemma 9.
\end{proof}

Now we are in a position to prove the main result of this section.

\vspace{3mm}
\noindent \textbf{Proposition 1.} {\it If $\widetilde{\mathfrak{A}} = \mathfrak{A}$, then $\mathfrak{A}$ possesses a bounded smooth strictly plurisubharmonic function.}

\begin{proof}
It follows from Lemma 9 that for each $n \in \mathbb N$ we can apply Theorem 4 to the set $K_n$. This way we will get a non-negative smooth plurisubharmonic function $\phi_n$ defined on the whole of ${\mathbb C}^2$ such that $K_n = \{p \in {\mathbb C}^2 : \phi_n(p) = 0 \}$ and $\phi_n$ is strictly plurisubharmonic on ${\mathbb C}^2 \setminus K_n$. Now if we define on the domain $\mathfrak{A}$ a new function
\[ \widetilde{\phi}_n(p) \coloneqq \left\{ \begin{array}{c@{\,,\quad}l} \phi_n(p)  &  {\rm{for}} \,\,  p \in {\widetilde{\mathfrak{A}}}_n, \\ 0 &   {\rm{for}} \,\,  p \in \mathfrak{A} \setminus {\widetilde{\mathfrak{A}}}_n,   \end{array} \right. \]
then, by construction of the set $K_n$, this function will still be a smooth bounded plurisubharmonic function on $\mathfrak{A}$ which is strictly plurisubharmonic exactly on the domain ${\widetilde{\mathfrak{A}}}_n$. Hence, for a decreasing sequence of positive numbers $\{ \varepsilon_n \}$ converging to zero fast enough, the function $\widetilde{\phi} \coloneqq \sum_{n=1}^\infty \varepsilon_n \widetilde{\phi}_n$ will be bounded smooth and plurisubharmonic on $\mathfrak{A}$ and, moreover, it will be strictly plurisubharmonic on the domain $\widetilde{\mathfrak{A}} = \bigcup_{n \in \mathbb{N}}{\widetilde{\mathfrak{A}}}_n$. Since, by our assumptions, $\widetilde{\mathfrak{A}} = \mathfrak{A}$, this completes
the proof of Proposition 1.
\end{proof}

\section{The case when $\widetilde{\mathfrak{A}} \neq \mathfrak{A}$}

\noindent
We first consider the following two sets 
\[ \mathcal E \coloneqq \mathfrak{A} \setminus \widetilde{\mathfrak{A}} \neq \emptyset \]
and
\[ {\mathcal E}^\ast \coloneqq p(\mathcal E), \]
where $p : {\mathbb C}_{z, w}^2 \to {\mathbb C}_z \times {\mathbb R}_u$ is the
canonical projection.

Then the special form of the domain $\mathfrak{A}$ and the construction of the domain $\widetilde{\mathfrak{A}}$ imply that:
\[ For \,\, every \,\, point \,\, (z, w) \in \mathcal E \,\, and \,\, every \,\, t > 0
\,\, the \,\, inclusion \,\, (z, w + it) \in \mathcal E \,\, holds.
 \]
Hence, we can define the set

\[ {\mathcal E}_0 \coloneqq \{(z, w) \in {\mathbb C}^2 : \,\, for \,\, every \,\, t
> 0 \,\, one \,\, has \,\,  (z, w + it) \in \mathcal E, \,\, but \,\, (z, w - it)
\notin \mathcal E\} \]
and then for each $(z, u) \in {\mathcal E}^\ast$ we denote by $v^\ast(z, u)$ the
real number such that $(z, u + iv^\ast(z, u)) \in {\mathcal E}_0$.

Let us also consider the set

\[ \mathcal F \coloneqq \{p \in {\mathbb C}^2 : p = \lim_{n \to \infty}q_n \,\, for
\,\, some \,\, sequence \,\, of \,\, points \,\, q_n \in \Gamma(F_n)\}, \]
where $\Gamma(F_n)$ is the graph of the function $F_n$ defined by (2) with $R = n$.

It follows directly from the definitions of these sets that the set $\mathcal E \cup
{\mathcal E}_0$ is closed and then from monotonicity of the domains
$\{\widetilde{\mathfrak{A}}_n\}_{n \in \mathbb N}$ (see Lemma 8) that $\mathcal F \subset
\mathcal E \cup {\mathcal E}_0$.

Now we describe how to pass to the limit for holomorphic leaves of the Levi
foliations of $\Gamma(\Phi_n)$ in order to get holomorphic disks in the limit set
${\mathcal F}$.

\vspace{0,5cm}
\textsc{Construction of the limit disk passing through the given point of $\mathcal F$.}

\vspace{0,2cm}
Let $(z_0, w_0)$ be an arbitrary point of the set ${\mathcal  F}$ and $(z_n, w_n) \in  \Gamma(F_n)$ be a sequence of points converging to $(z_0, w_0)$. For each $n = 0, 1, 2, ...$  we consider the cylinder

\[ {\mathcal C}_n \coloneqq \{(z, u) : |z - z_n| < 1, u_n - 1 < u < u_n + 1 \} \subset {\mathbb C}_z \times {\mathbb R}_u, \]
where $w_n = u_n + iv_n$. Then, for $n$ sufficiently large (so large that $ {\mathcal C}_n \subset B_n(0)$ and, hence, the function $\Phi_n$ is defined on ${\mathcal C}_n$), we denote by ${\mathcal L}^\alpha_n$ the leaves of the Levi foliation of the set $\Gamma(\Phi_n) \cap ({\mathcal C}_n \times {\mathbb R}_v)$. If for each leaf ${\mathcal L}^\alpha_n$ of this foliation we consider its projection $p({\mathcal L}^\alpha_n)$ to the cylinder $ {\mathcal C}_n$ (here $p$ is the canonical projection $p : {\mathbb C}^2_{z, w} \to {\mathbb C}_z \times {\mathbb R}_u$) and apply to this surface Lemma 2, we will obtain a domain ${\Omega}^\alpha_n \subset {\mathbb C}_z$ and a holomorphic function $h^\alpha_n : {\Omega}^\alpha_n \to {\mathbb C}_w$ such that its graph $\Gamma(h^\alpha_n)$ coincides with the leaf ${\mathcal L}^\alpha_n$. Consider now for each $n$ the special leaf ${\mathcal L}^{\alpha_0}_n$ which has the property $p({\mathcal L}^{\alpha_0}_n) \ni (z_n, u_n)$ and observe that, by Lemma 3 (applied to ${\Omega}^{\alpha_0}_n$ on the place of $\Omega$, $U^{\alpha_0}_n \coloneqq \Re h^{\alpha_0}_n$ on the place of $U$, $V^{\alpha_0}_n \coloneqq \Im h^{\alpha_0}_n$ on the place of $V$ and $ {\mathcal C}_n$ on the place of $G$), we have that the domain ${\Omega}^{\alpha_0}_n \ni z_n$ is simply-connected. It follows now from Lemma 4 (applied to ${\Omega}^{\alpha_0}_n$ on the place of $\Omega$, $U^{\alpha_0}_n$ on the place of $U$, $V^{\alpha_0}_n$ on the place of $V$ and with $M_n \coloneqq \inf_{(z, u) \in {\mathcal C}_n}F(z, u)$ and $L_n \coloneqq v_n$), that there exists $R^\ast_n > 0$ which only depends on $M_n$ and $L_n$ and such that the disk $\Delta_{R^\ast_n}(z_n) \coloneqq  \{ z \in {\mathbb C} : |z - z_n| < {R^\ast_n} \}$ is contained in ${\Omega}^{\alpha_0}_n$. Moreover, the arguments used in the proof of Lemma 4 (see for details the proof of Lemma 3.5 in \cite{Shcherbina93}) show that $R^\ast$ there depends on $M$ and $L$ in a continuous way. Then, since by our assumptions $\lim_{n \to \infty}v_n = v_0$, and since by continuity of $F$ one also has that $\lim_{n \to \infty}M_n = M_0$, we conclude that for each $R < R^\ast_0 \eqqcolon \lim_{n \to \infty}R^\ast_n$ the sequence of holomorphic functions $h^{\alpha_0}_n$ is defined on the disk $\Delta_R(z_0)$ for all $n$ large enough. Moreover, since for each $n$ we have that $|U^{\alpha_0}_n(z) - u_n| < 1$ for $z \in \Delta_{R^\ast_n}(z_0)$, and since $M_0 \leftarrow M_n < V^{\alpha_0}_n$ on $\Delta_{R^\ast_n}(z_0)$, and $V^{\alpha_0}_n(z_n) \to v_0$ as $n \to \infty$, we can apply Montel theorem and choose a subsequence $h^{\alpha_0}_{n_k}$ which converges to a holomorphic function $h_0$ uniformly on compact subsets of $\Delta_{R^\ast_0}(z_0)$. Finally, we observe that by construction one has that $\Im h_0(z_0) = v_0$. Then the graph $\Gamma(h_0)$ of the function $h_0$ is the desired holomorphic disk which will be denoted in what follows by $D_0$. 

\vspace{3mm}
In the next statement we collect the properties of the described above limit disk which we will need later on.

\vspace{3mm}
\noindent \textbf{Lemma 10.} {\it For each point $(z_0, w_0) \in {\mathcal  F}$ there is a holomorphic disk $D_0$ of the form
\[ D_0 = \{z, w) \in {\mathbb C}^2 : w = h_0(z), z \in \Delta_{R^\ast_0}(z_0) \},
 \]
obtained by the described above limit procedure and such that $(z_0, w_0) \in D_0 \subset {\mathcal F}$ with $R^\ast_0$ which depends only on the estimate from above on $v_0$ and on the estimate from below on $\inf_{(z, u) \in {\mathcal C}_0}F(z, u)$. 

Moreover, the following properties of the constructed limit disk hold:
\begin{enumerate}
  \item[(A)]  The disk $D_0$ is uniquely determined (even locally). 

  \item[(B)]  If $(z_0, w_0) \in {\mathcal E}_0$, then $D_0 \subset {\mathcal E}_0$.
\end{enumerate}

The property (A) means, more precisely, that: For each $r \leq R^\ast_0$ and each holomorphic disk ${\widetilde D}_0 \ni (z_0, w_0)$ of the form
\[ {\widetilde D}_0 = \{z, w) \in {\mathbb C}^2 : w = {\widetilde h}_0(z), z \in \Delta_r(z_0) \},
 \]
which is also obtained by the limit procedure above only applied to another subsequence $h^{\alpha_0}_{n_m}$ of the sequence $h^{\alpha_0}_n$, which converges to a holomorphic function ${\widetilde h}_0$ uniformly on compact subsets of $\Delta_r(z_0)$, we have ${\widetilde D}_0 \subset D_0$.}

\begin{proof}
To prove Claim (A) we argue by contradiction and assume that there is a holomorphic disk ${\widetilde D}_0$ of the form described in Lemma 10 such that  ${\widetilde D}_0 \not \subset D_0$. Then, since $D_0 \cap {\widetilde D}_0 \ni (z_0, w_0)$, we conclide from Rouche's theorem, that for $k$ and $m$ large enough and $t > 0$ small enough the disks 

\[ D^{\alpha_0}_{n_k} \coloneqq \Gamma(h^{\alpha_0}_{n_k}|\Delta_r(z_0)) \subset \Gamma(\Phi_{n_k}) \]
and 
\[ {\widetilde D}^{\alpha_0}_{{n_m},t} \coloneqq \Gamma((h^{\alpha_0}_{n_m} + it)|\Delta_r(z_0)) \subset \Gamma(\Phi_{n_m} + t), \]
which approximate the part of the disks $D_0$ over $\Delta_r(z_0)$ and the disk ${\widetilde D}_0$, respectively, will also have a non-empty intersection:

\[ D^{\alpha_0}_{n_k} \cap {\widetilde D}^{\alpha_0}_{{n_m},t} \neq \emptyset. \]
On the other hand, if $n_m > n_k$, then it follow from Lemma 7 that for $t > 0$ one has that $\Phi_{n_k} < \Phi_{n_m} + t$ on $B_{n_k}(0)$, and hence also that $\Gamma(\Phi_{n_k}) \cap \Gamma(\Phi_{n_m} + t) = \emptyset$. We conclude now that 

\[ D^{\alpha_0}_{n_k} \cap {\widetilde D}^{\alpha_0}_{{n_m},t} =   \Gamma(h^{\alpha_0}_{n_k}|\Delta_r(z_0)) \cap \Gamma((h^{\alpha_0}_{n_m} + t)|\Delta_r(z_0)) \subset \Gamma(\Phi_{n_k}) \cap \Gamma(\Phi_{n_m} + t) = \emptyset, \]
which gives us a desired contradiction.

To prove Claim (B) we first observe that, in view of the definitions of $D_0$ and
$\mathcal F$, we obviously have the inclusion $D_0 \subset \mathcal F$. Then, since
we also have the inclusion $\mathcal F \subset \mathcal E \cup {\mathcal E}_0$, we
conclude that for each $z \in \Delta_{R^\ast_0}(z_0)$ one has that $(z, \Re h_0(z)) \in
{\mathcal E}^\ast$ and $v^\ast(z,\Re h_0(z)) \leq \Im h_0(z)$. Hence, for proving
that $v^\ast(z,  \Re h_0(z)) = \Im h_0(z)$ (which gives us the desired statement that $D_0
\subset {\mathcal E}_0$), it is enough to show that for each $z \in
\Delta_{R^\ast_0}(z_0)$ one has that $v^\ast(z, \Re h_0(z)) \geq \Im h_0(z)$. We argue by
contradiction and assume that for some $z \in \Delta_{R^\ast_0}(z_0)$ we have that
$v^\ast(z, \Re h_0(z)) < \Im h_0(z)$. Then, if we define 

\[ r^\ast \coloneqq \sup \{ r \in [0, R^\ast_0) : v^\ast(z, \Re h_0(z)) = \Im h_0(z) {\rm \,\, for \,\, all \,\,} z \in \overline{\Delta}_r(z_0) \}, \]
we see that $r^\ast < R^\ast_0$ and that we can find two arbitrary close to each other
points $z' \in \overline{\Delta}_{r^\ast}(z_0)$ and $z'' \notin \overline{\Delta}_{r^\ast}(z_0)$ such that
one has $v^\ast(z', u') = \Im h_0(z') \eqqcolon v'$ and $v^\ast(z'', u'') < \Im
h_0(z'') \eqqcolon v''$, where $u' \coloneqq \Re h_0(z')$ and $u'' \coloneqq \Re
h_0(z'')$, respectively.

Let us now take some number $\tilde r \in (r^\ast, R^\ast_0)$ which we will keep fixed till the end of the proof of Lemma 10. Since $D_0 \subset \mathcal F$, we have the following inclusion 
\[ E_{\tilde{r}} \coloneqq \Gamma(\Re h_0 | \overline{\Delta}_{\tilde{r}}(z_0))  = p(\Gamma(h_0) | \overline{\Delta}_{\tilde{r}}(z_0)) = p(D_0 \cap (\overline{\Delta}_{\tilde{r}}(z_0) \times {\mathbb C}_w)) \subset {\mathcal E}^\ast. \]
Then we define 
\[ L_{\tilde r} \coloneqq \max \{\Im h_0(z) : z \in \overline{\Delta}_{\tilde{r}}(z_0) \}, \]
and 
\[ M_{\tilde r} \coloneqq \inf_{(\tilde{z}, \tilde{u}) \in E_{\tilde{r}}}\inf_{(z, u) \in {{\mathcal C}(\tilde{z}, \tilde{u})}}F(z, u), \]
where 
\[ {\mathcal C}(\tilde{z}, \tilde{u}) \coloneqq \{(z, u) : |z - \tilde{z}| < 1, \tilde{u} - 1 < u < \tilde{u} + 1 \} \subset {\mathbb C}_z \times {\mathbb R}_u, \]
and observe that, by continuity of $h_0$ and $F$, one has $- \infty <  M_{\tilde r} \leq L_{\tilde r} < + \infty$.

It follows then from our construction above of the limit disk passing  through the given point of $\mathcal F$, and in view of the uniform estimates on $v_0$ from above by $L_{\tilde r}$ and on $\inf_{(z, w) \in {\mathcal C}_0}F(z, w)$ from below by $M_{\tilde r}$ for points 
\[ (z_0, w_0) \in {\mathcal F}_{\tilde r} \coloneqq \{(z, w) : (z, u)  \in E_{\tilde{r}}, v^\ast(z, u) \leq v \leq \Im h_0(z) \}, \]
that we also have the inequality
\[ \widetilde R \coloneqq \inf_{(z_0, w_0) \in {\mathcal F}_{\tilde r}} \{R^\ast_0(z_0, w_0) \} > 0, \]
where $R^\ast_0(z_0, w_0)$ denote the radius $R^\ast_0$ of our construction above corresponding to the point $(z_0, w_0)$.

Now we can finally specify our choice of the described above points $z' \in \overline{\Delta}_{r^\ast}(z_0)$ and $z'' \notin \overline{\Delta}_{r^\ast}(z_0)$. Namely, we first define 
\[ \tilde{r}^\ast \coloneqq \frac{1}{3}\min\{\widetilde R, \tilde{r} - r^\ast\} \]
and then choose $z'$ and $z''$ so close to each other that $|z' - z''| < \tilde{r}^\ast$. 

Consider now the limit disks $D'$ and $D''$ that correspond by our construction above to the points $(z', \Re h_0(z') + iv^\ast(z', \Re h_0(z'))$ and $(z'', \Re h_0(z'') + iv^\ast(z'', \Re h_0(z''))$, respectively. By our choice of $\tilde{r}^\ast$, after maybe shrinking the disks $D'$ and $D''$ if necessary, we can assome that the disk $D'$ is the graph of a holomorphic function $w = h'(z)$ defined on the disk $\Delta_{2\tilde{r}^\ast}(z')$ and the disk $D''$ is the graph of a holomorphic function $w = h''(z)$ defined on the disk $\Delta_{2\tilde{r}^\ast}(z'')$. Observe that, by our choice of the point $z'$, and in view of the unicity property of the limit disks established in Part (A) of this lemma, we have that $h'(z) = h_0(z)$ for all $z \in \Delta_{2\tilde{r}^\ast}(z')$. We consider now the following two cases:
\begin{enumerate}
  \item[(i)] $\Re h_0(z) = \Re h''(z) {\rm \,\, for \,\, all \,\, } z \in \Delta_{2\tilde{r}^\ast}(z') \cap \Delta_{2\tilde{r}^\ast}(z''),$

  \item[(ii)] $\Re h_0(z^\ast) \neq \Re h''(z^\ast) {\rm \,\, for \,\, some \,\, } z^\ast \in \Delta_{2\tilde{r}^\ast}(z') \cap \Delta_{2\tilde{r}^\ast}(z'').$
\end{enumerate}

In the Case (i) we have that $\Im h''(z) = \Im h_0(z) + C$ for some constant $C$ and all $z \in \Delta_{2\tilde{r}^\ast}(z') \cap \Delta_{2\tilde{r}^\ast}(z'')$. Since, by our choice of the point $z''$, we have that $v^\ast(z'', u'') < \Im h_0(z'')$, it follows that $C < 0$. But then, by the choice of the point $z'$,  
\[ \Im h''(z') = \Im h_0(z') + C < \Im h_0(z') = v^\ast(z', \Re h_0(z')). \]
On the other hand, since $D'' \subset \mathcal F$, we conclude from the definition of the function $v^\ast(z, u)$ that
\[ v^\ast(z', \Re h_0(z')) \leq  \Im h''(z'). \]
The last two inequalities contradict each other. This gives us a desired contradiction in the Case (i).

In the Case (ii) we use an argument similar to the one we used in the proof of Part (A) of this lemma. Namely, by our construction above of the limit disks passing through the given point of $\mathcal F$, and in view of the unicity of the limit disks, we can consider two sequences of holomorphic disks $D_{0, n}$ and $D'_m$ of the form,

\begin{equation} \label{D_o,n}
D_{0, n} \coloneqq \Gamma(h_{0, n}) \subset \Gamma(\Phi_n)
\end{equation}
and 

\begin{equation} \label{D''_m}
D''_m \coloneqq \Gamma(h''_m) \subset \Gamma(\Phi_m),
\end{equation}
which approximate the part of the disks $D_0$ over $\Delta_{2\tilde{r}^\ast}(z'')$ and the disk $D''$, respectively. The functions $h_{0, n}$ and $h''_m$ here are defined and holomorphic on the disk $\Delta_{2\tilde{r}^\ast}(z'')$ and have the property that $\lim_{n \to \infty}h_{0, n}(z'') = h_0(z'')$ and $\lim_{m \to \infty}h''_m(z'') = h''(z'')$. 

Observe now that by our choice of the point $z''$ and the disk $D''$ one has that 

\[ \Im h''(z'') = v^\ast(z'', \Re h_0(z'')) < \Im h_0(z''), \]
and then, in view of the uniform convergence on the disk $\Delta_{2\tilde{r}^\ast}(z'')$ of the sequences $h_{0, n}(z)$ and $h''_m(z)$ to $h_0(z)$ and $h''(z)$, respectively, we conclude that

\begin{equation} \label{lim_ineq} 
\Im h''_m(z_l)  < \Im h_{0, n}(z_l)
\end{equation}
for every sequence $\{z_l\}_{n = 1}^\infty$ of points in $\Delta_{2\tilde{r}^\ast}(z'')$ converging to $z''$ and all $m$, $n$ and $l$ large enough.

On the other hand, since, by the choice of the disk $D''$, we have $\Re h''(z'') = \Re h_0(z'')$ and since, by the conditions in the Case (ii), we have $\Re h'' \neq \Re h_0$ in $\Delta_{2\tilde{r}^\ast}(z') \cap \Delta_{2\tilde{r}^\ast}(z'')$ we conclude that for the functions $\Re h''_m$ and $\Re h_{0, n}$ approximating the functions $\Re h''$ and $\Re h_0$, respectively, we have in the case when $n$ and $m$ are large enough that there exists a point $z_{n, m} \in \Delta_{2\tilde{r}^\ast}(z') \cap \Delta_{2\tilde{r}^\ast}(z'')$ as close to $z''$ as we wish and such that  $\Re h''_m(z_{n, m}) = \Re h_{0, n}(z_{n, m})$. Hence, if in the last statement for each $n$ we in addition choose the value of $m$ in such a way that $m > n$, then, by Lemma 7, and, in the view of inclusions (5) and (6), we will have that 

\[ \Im h''_{m}(z_{n, m}) = \Phi_{m}(z_{n, m}, \Re h''_{m}(z_{n, m})) \geq \Phi_n(z_{n, m}, \Re h_{0, n}(z_{n, m})) = \Im h_{0, n}(z_{n, m}). \]
The last inequality contradicts the inequality (7). This completes the proof of Lemma 10.
\end{proof}

Now let us recall some definitions.

\vspace{3mm}
Let $\Sigma$ be a domain over $\mathbb C$ with projection $\pi : \Sigma \to \mathbb C$, and let $f = U + iV$ be a smooth function in $\Sigma$.

\vspace{0,3cm}
\noindent \textbf{Definition 4.} Let $d$ be the pull-back on $\Sigma$ of the standard metric on $\mathbb C$. Define the distance $\rho(\xi_1, \xi_2)$ between two arbitrary points $\xi_1$ and $\xi_2$ in $\Sigma$ as $\rho(\xi_1, \xi_2) \coloneqq \inf_\gamma{l(\gamma)}$, where $l(\gamma)$ is the length of the curve $\gamma$ with respect to the metric $d$, and the infimum is taken over all smooth curves in $\Sigma$ connecting the points $\xi_1$ and $\xi_2$. 

\vspace{0,3cm}
\noindent \textbf{Definition 5.} Let $G$ be a domain of $\mathbb C \times \mathbb R$. The graph $\Gamma(U) = \{(\pi(\xi), U(\xi)) \in \mathbb C \times \mathbb R : \xi \in \Sigma \}$ is said to be {\it locally embedded} in
$G$ if $\Gamma(U) \subset G$, and for each $\xi_1, \xi_2 \in \Sigma$, such that $\xi_1 \neq \xi_2$, one has $(\pi(\xi_1), U(\xi_1)) \neq (\pi(\xi_2), U(\xi_2))$.

\vspace{0,3cm}
\noindent \textbf{Remark 1.} The property of $\Gamma(U)$ to be locally embedded in $G$ differs from the usual embedding property. For example, a surface $\Gamma(U)$, locally embedded according to our definition, might be an everywhere dense subset of $G$ and this, obviously, cannot happen for embedded surfaces.

For each point $\xi \in \Sigma$, denote by $\rho(\xi, \partial{\Sigma})$ the maximum of all real numbers $\rho$ for which there exists a continuous map $U'$ from the disc $\Delta_\rho(\pi(\xi)) = \{z \in \mathbb C: |z - \pi(\xi)| < \rho \}$ to $\Sigma$ with the following properties: $\pi \circ U'$ is the identity map on $\Delta_\rho(\pi(\xi))$ and $U'(\pi(\xi)) = \xi$.

\vspace{0,3cm}
\noindent \textbf{Definition 6.} We say that the sequence $\{\xi_n \}$ of points of $\Sigma$ {\it converges to the boundary} $\partial{\Sigma}$ as $n \to + \infty$ (and then we write $\xi_n \to \partial{\Sigma}$) if $\rho(\xi, \partial{\Sigma}) \to 0$ → 0 as $n \to + \infty$.

\vspace{0,3cm}
\noindent \textbf{Definition 7.} We define the boundary $\partial{\Gamma}(U)$ of the surface ${\Gamma}(U)$ as the cluster set of all sequences $\{(\pi(\xi_n), U(\xi_n))\}$ such that $\xi_n \to \partial{\Sigma}$ as $n \to + \infty$.

\vspace{0,3cm}
\noindent \textbf{Definition 8.} Let $G$ be a domain of $\mathbb C \times \mathbb R$. We say that the surface ${\Gamma}(U)$ is {\it properly embedded} into the domain $G$ if $\partial{\Gamma}(U) \cap G = \emptyset$.

Similarly, we define the boundary $\partial{\Gamma}(f)$ of the surface ${\Gamma}(f) \coloneqq \{(\pi(\xi), f(\xi)) \in {\mathbb C}^2 : \xi \in \Sigma \}$ and the property of $\Gamma(f)$ to be properly embedded into a domain $\mathcal G \subset {\mathbb C}^2$.

\vspace{0,3cm}
Now we describe how to extend the limit disk $D_0$ from the Part (B) of Lemma 10 to the maximal limit leaf which is contained in the set ${\mathcal E}_0$.

\vspace{0,5cm}
\textsc{Construction of the maximal limit leaf ${\mathcal L}_0$ contained in  ${\mathcal E}_0$.}

\vspace{0,2cm}
Let $(z_0, w_0)$ be an arbitrary point of the set ${\mathcal  E}_0$. We define ${\mathcal L}_0$ as the set of all points $(z', w') \in {\mathcal  E}_0$ such that there are finitely many points $(z_1, w_1), (z_2, w_2), ..., (z_m, w_m) \in {\mathcal  E}_0$, with $(z_m, w_m) = (z', w')$, satisfying the following property:

\vspace{2mm}

\vspace{2mm}
\noindent
{\it For each $i = 0, 1, 2, ..., m - 1$, one has that $(z_{i + 1}, w_{i + 1}) \in D_i$, where $D_i$ is the disk which corresponds to the disk $D_0$ of Lemma 10, when we take the point $(z_i, w_i)$ on the place of $(z_0, w_0)$.}

Observe first that, since, by definition of ${\mathcal L}_0$, each point $(z', w') \in {\mathcal L}_0$ is contained in the corresponding disk $D_{m - 1} \subset {\mathcal  E}_0$, it follows that ${\mathcal L}_0$ with the projection to ${\mathbb C}_z$, defined as the restriction  to ${\mathcal L}_0$ of the canonical projection $p_w : {\mathbb C}_{z, w}^2 \to {\mathbb C}_z$, is a domain over ${\mathbb C}_z$. Then, from the unicity of the limit disks (see Part (A) of Lemma 10) we conclude that ${\mathcal L}_0$ is locally embedded into ${\mathbb C}_{z, w}^2$. Moreover, since, in view of the definitions of ${\mathcal  E}_0$ and ${\mathcal  E}^\ast$, the map $p : {\mathcal  E}_0 \to {\mathcal  E}^\ast$ is bijective, one has that the projection ${\mathcal M}_0 \coloneqq p({\mathcal L}_0)$ of ${\mathcal L}_0$ to ${\mathbb C}_z \times {\mathbb R}_u$  is locally embedded into ${\mathbb C}_z \times {\mathbb R}_u$. 

Observe further that ${\mathcal L}_0$ is properly embedded into ${\mathbb C}_{z, w}^2$. Indeed, if not, then there is a cluster point $(z^\ast, w^\ast)$ for the set ${\mathcal L}_0$, i.e. there exists a sequence of points $(z_n, w_n)$ of ${\mathcal L}_0$ (as a domain over ${\mathbb C}_z$) such that $(z_n, w_n) \to (z^\ast, w^\ast)$ as $n \to \infty$ and
\begin{equation} \label{d_L_0} 
\lim_{n \to + \infty}\rho((z_n, w_n), \partial{{\mathcal L}_0}) \to 0. 
\end{equation}
On the other hand, since $v^\ast \coloneqq \Im w^\ast < + \infty$, and since $v_n \coloneqq \Im w_n \to \Im w^\ast$ as $n \to + \infty$, we conclude from Lemma 10 that 
\begin{equation} \label{R_ast} 
\liminf_{n \to \infty} R_n^\ast \geq R^\ast > 0,
\end{equation}
where $R_n^\ast$ and $R^\ast$ represent the radius $R_0^\ast$ of Lemma 10 computed at the points $(z_n, w_n)$ and $(z^\ast, w^\ast)$, respectively. Since for each $n \in \mathbb N$ we obviously have that
\[ \rho((z_n, w_n), \partial{{\mathcal L}_0}) \geq R_n^\ast, \]
it follows from (9) that 
\[ \liminf_{n \to \infty}\rho((z_n, w_n), \partial{{\mathcal L}_0}) \geq R^\ast > 0. \]
The last inequality contadicts to (8). This proves that ${\mathcal L}_0$ is properly embedded into ${\mathbb C}_{z, w}^2$.

\vspace{3mm}
In the next lemma we prove that the embedding of ${\mathcal L}_0$ in ${\mathbb C}_{z, w}^2$ is actually much more special.

\vspace{3mm}
\noindent \textbf{Lemma 11.} {\it The surface ${\mathcal L}_0$ is properly embedded into ${\mathbb C}^2$ and projects one-to-one on its image in ${\mathbb C}_z$. More precisely, there is a domain $\mathfrak{D} \subset {\mathbb C}_z$ and a holomorphic function ${\mathcal H} : \mathfrak{D} \to {\mathbb C}_w$ such that the graph $\Gamma({\mathcal H})$ of ${\mathcal H}$ coincides with ${\mathcal L}_0$ and, moreover, $|{\mathcal H}(z)| \to + \infty$ as $z \to \partial{\mathfrak{D}}$.} 

\begin{proof}
First, we prove that the surface ${\mathcal L}_0$ projects one-to-one on its image in ${\mathbb C}_z$. For doing this we argue by contradiction and assume that this property does not hold. Since ${\mathcal L}_0$ projects one-to-one to ${\mathcal M}_0 = p({\mathcal L}_0) \subset {\mathbb C}_z \times {\mathbb R}_u$, our assumption implies that there are two points $q' = (\tilde{z}, \tilde{u}'), q'' = (\tilde{z}, \tilde{u}'') \in {\mathcal M}_0$ with $\tilde{u}' \neq \tilde{u}''$. Since, by construction, the surface ${\mathcal L}_0$ is connected, we conclude that ${\mathcal M}_0$ is also connected, and hence there is a smooth curve $\gamma \subset {\mathcal M}_0$ connecting points $q'$ and $q''$. Moreover, since $p^{-1}(\gamma) \subset {\mathcal L}_0$ is compact, for points $(z^\ast, w^\ast) \in p^{-1}(\gamma)$, $\, w^\ast = u^\ast + iv^\ast,$ taken on the place of $(z_0, w_0)$ in Lemma 10 one has uniform bounds from above on $v^\ast$ and from below on $\inf_{(z, w) \in {\mathcal C}^\ast}F(z, w)$, where 

\[ {\mathcal C}^\ast \coloneqq \{(z, u) : |z - z^\ast| < 1, u^\ast - 1 < u < u^\ast + 1 \} \subset {\mathbb C}_z \times {\mathbb R}_u. \]
It follows then from Lemma 10 that there is a number $R(\gamma) > 0$ such that at each point $(z^\ast, w^\ast) \in p^{-1}(\gamma)$ for the radius $R^\ast_0(z^\ast, w^\ast)$, which is the radius $R^\ast_0$ of Lemma 10, but computed at the point $(z^\ast, w^\ast)$ on the place of the point $(z_0, w_0)$, we have the inequality $R^\ast_0(z^\ast, w^\ast) \geq R(\gamma)$. 

Let us now take a finite set of points $Q_0, Q_1, Q_2, ..., Q_m \in p^{-1}(\gamma)$ such that:
\begin{enumerate}
  \item[(1)] $\,$ $p(Q_0) = q'$ and $p(Q_m) = q''$.

  \item[(2)] $\,$ For each $i = 0, 1, 2, ..., m - 1$ the inclusion $p_u(\gamma_i) \subset \Delta_{R(\gamma)}(z_i) \cap \Delta_{R(\gamma)}(z_{i + 1})$ holds.
\end{enumerate}

\noindent
Here, for each $i = 0, 1, ..., m$, $Q_i \eqqcolon  (z_i, w_i)$, for each $i = 0, 1, ..., m - 1$, $\gamma_i$ is the part of $\gamma$ between $p(Q_i)$ and $p(Q_{i + 1})$, and $p_u : {\mathbb C}_z \times {\mathbb R}_u \to {\mathbb C}_z$ is the canonical projection.

\vspace{0,3cm}
It follows then from our choice of $R(\gamma)$ that for each $n$ sufficiently large there exists a holomorphic function $h^0_n : \Delta_{R(\gamma)}(z_0) \to {\mathbb C}_w$ such that $\Re h^0_n(z_0) = \Re w_0$ and the graph $\Gamma(h^0_n)$ of the function $h^0_n$ is contained in a holomorphic leaf of the Levi foliation on $\Gamma(\Phi_n)$. Moreover, in view of Lemma 10, the functions $h^0_n$ converge as $n \to \infty$ uniformly on compact subsets of $\Delta_{R(\gamma)}(z_0)$ to a holomorphic function $h^0_0 : \Delta_{R(\gamma)}(z_0) \to {\mathbb C}_w$ such that $h^0_0(z_0) = w_0$ and $\Gamma(h^0_0) \subset {\mathcal L}_0$. Then, by the Condition (2) above, and by the choice of the point $Q_1$, we also see that $h^0_0(z_1) = w_1$.

Since $Q_1 = (z_1, w_1) \in p^{-1}(\gamma)$, and since $h^0_n(z_1) \to w_1$ as $n \to \infty$, it follows once again from our choice of $R(\gamma)$ that for each $n$ sufficiently large there exists a holomorphic function $h^1_n : \Delta_{R(\gamma)}(z_1) \to {\mathbb C}_w$ such that $\Re h^1_n(z_1) = \Re h^0_n(z_1)$ and the graph $\Gamma(h^1_n)$ of the function $h^1_n$ is contained in a holomorphic leaf of the Levi foliation on $\Gamma(\Phi_n)$. Observe, that the last conditions imply that $h^1_n = h^0_n$ on the intersection $\Delta_{R(\gamma)}(z_0) \cap \Delta_{R(\gamma)}(z_1)$ of the disks $\Delta_{R(\gamma)}(z_0)$ and $\Delta_{R(\gamma)}(z_1)$. Then, as before, in view of Lemma 10 again, the functions $h^1_n$ converge as $n \to \infty$ uniformly on compact subsets of $\Delta_{R(\gamma)}(z_1)$ to a holomorphic function $h^1_0 : \Delta_{R(\gamma)}(z_1) \to {\mathbb C}_w$ such that $h^1_0(z_1) = w_1$ and $\Gamma(h^1_0) \subset {\mathcal L}_0$. Hence, by the unicity property (A) of Lemma 10, we see that $h^1_0 = h^0_0$ on the set $\Delta_{R(\gamma)}(z_0) \cap \Delta_{R(\gamma)}(z_1)$ and, by the Condition (2) above, and by the choice of the point $Q_2$, we also see that $h^1_0(z_2) = w_2$.

If we repeat inductively this argument for each $i = 2, 3, ..., m - 1$, we will get for all $n$ large enough holomorphic functions $h^i_n : \Delta_{R(\gamma)}(z_i) \to {\mathbb C}_w$ such that $\Re h^i_n(z_i) = \Re h^{i - 1}_n(z_i)$ and the graph $\Gamma(h^i_n)$ of each function $h^i_n$ is contained in a holomorphic leaf of the Levi foliation on $\Gamma(\Phi_n)$. The last conditions actually imply that $h^i_n = h^{i - 1}_n$ on the intersection $\Delta_{R(\gamma)}(z_{i - 1}) \cap \Delta_{R(\gamma)}(z_i)$ of the disks $\Delta_{R(\gamma)}(z_{i - 1})$ and $\Delta_{R(\gamma)}(z_i)$. Then, in view of Lemma 10 again, each sequence of functions $h^i_n$ will  converge as $n \to \infty$ uniformly on compact subsets of $\Delta_{R(\gamma)}(z_i)$ to a holomorphic function $h^i_0 : \Delta_{R(\gamma)}(z_i) \to {\mathbb C}_w$ such that $h^i_0(z_i) = w_i$ and $\Gamma(h^i_0) \subset {\mathcal L}_0$. Hence, by the unicity property (A) of Lemma 10, we see that $h^i_0 = h^{i - 1}_0$ on the set $\Delta_{R(\gamma)}(z_i) \cap \Delta_{R(\gamma)}(z_{i - 1})$ and, by the  Condition (2) above, and by the choice of the point $Q_{i + 1}$, we  see that $h^i_0(z_{i + 1}) = w_{i + 1}$.

Let us now denote 

\begin{equation} \label{h_m}
\delta \coloneqq |\tilde{u}' - \tilde{u}''| > 0. 
\end{equation}
Then, by uniform convergence of the functions $h^{m - 1}_n$ to $h^{m - 1}_0$ on compact subsets of $\Delta_{R(\gamma)}(z_{m - 1})$ as $n \to \infty$, and in view of the inclusion $z_{m} \in \Delta_{R(\gamma)}(z_{m - 1})$ (see Condition (2) above with $i = m - 1$), we can take $n$ so large that 
\begin{equation} \label{h__m} 
|\Re h^{m - 1}_n(z_m) - \Re h^{m - 1}_0(z_m)| < \frac{\delta}{2}.
\end{equation}
Observe, that from our construction above one has

\[ p^{- 1}(\gamma) \subset \bigcup_{1 \leq i \leq m - 1}{\Gamma(h^i_0)} \subset  {\mathcal L}_0, \]
and hence, from the choice of the points $Q_0$, $ Q_m$ and points $z_0$, $z_m$, we can conclude that $Q_0 = (z_0, h^0_0(z_0))$ and $Q_m = (z_m, h^{m - 1}_0(z_m))$. It follows then from the choice of the points $q'$, $q''$ and the Condition (1) above that $\Re h^{m - 1}_0(z_m) - \Re h^0_0(z_0) = \tilde{u}'' - \tilde{u}'$ and, therefore, by (10),

\begin{equation} 
|\Re h^{m - 1}_0(z_m) - \Re h^0_0(z_0)| = \delta.
\end{equation}

\noindent
By the choice of the points $Q_0$, $ Q_m$ and points $z_0$, $z_m$, and in view of the Condition (1), we can also see that $z_0 = \tilde{z}  = z_m$. Then, we can finally obtain from (11) and (12) that

\begin{equation} \begin{split}
|\Re h^{m - 1}_n(z_m) - \Re h^0_n(z_0)| = |\Re h^{m - 1}_n(z_m) - \Re h^0_0(z_0)| &\geq \\ |\Re h^{m - 1}_0(z_m) - \Re h^0_0(z_0)| - |\Re h^{m - 1}_n(z_m) - \Re h^{m - 1}_0(z_m)| &> \frac{\delta}{2}.
\end{split} \end{equation}
\noindent
Denote now by ${\mathcal L}_n$ the holomorphic leaf of the Levi foliation on $\Gamma(\Phi_n)$ such that $p({\mathcal L}_n) \ni q'$. Then observe that the points ${\tilde{q}}' \coloneqq q' = (z_0, \Re h^0_n(z_0))$ and ${\tilde{q}}'' \coloneqq (z_m, \Re h^{m - 1}_n(z_m))$ belong to the surface ${\mathcal M}_n \coloneqq p({\mathcal L}_n)$ and that $z_0 = z_m$, but according to (13) we have $\Re h^0_n(z_0)) \neq \Re h^{m - 1}_n(z_m))$.  This means that on the surface ${\mathcal M}_n$, which is smooth, properly embedded into the ball $B_n(0)$ and locally has one-to-one projection on its image in ${\mathbb C}_z$, we have two different points that project to the same point of ${\mathbb C}_z$. According to Lemma 2 above this is not possible. Hence, we arrived to a contradiction which proves our claim that the surface ${\mathcal L}_0$ projects one-to-one on its image in ${\mathbb C}_z$.

Let us now denote by $\mathfrak{D}$ the projection of ${\mathcal L}_0$ to ${\mathbb C}_z$ and then by ${\mathcal H} : \mathfrak{D} \to {\mathbb C}_w$ a holomorphic function whose graph $\Gamma({\mathcal H})$ coincides with the surface ${\mathcal L}_0$. Then the last statement of Lemma 11 claiming that $|{\mathcal H}(z)| \to + \infty$ as $z \to \partial{\mathfrak{D}}$ is just a reformulation of the properness of the embedding of ${\mathcal L}_0$ into ${\mathbb C}_{z, w}^2$ proved after the construction of the maximal leaf ${\mathcal L}_0$ just before the statement of Lemma 11. The proof of Lemma 11 is now completed.
\end{proof}

The next statement shows that the shape of the surface ${\mathcal L}_0$ is even more special than it is claimed in the last lemma.

\vspace{3mm}
\noindent \textbf{Lemma 12.} {\it The domain $\mathfrak{D}$ coincides with the whole of ${\mathbb C}_z$ and, therefore, ${\mathcal H} : {\mathbb C}_z \to {\mathbb C}_w$ is an entire function such that $\Gamma({\mathcal H}) = {\mathcal L}_0$. The last equality implies, in particular, that the inclusion $\Gamma({\mathcal H}) \subset \overline{\mathfrak{A}}$ holds.}

\begin{proof}
We first prove that the domain $\mathfrak{D}$ is simply-connected. Consider, as above, the surface ${\mathcal M}_0 = p({\mathcal L}_0) \subset {\mathbb C}_z \times {\mathbb R}_u$ and  for each sufficiently large $n \in \mathbb N$ let ${\mathcal M}_0^n$ be a connected component of the surface ${\mathcal M}_0 \cap B_n(0)$ containing a fixed point, say the point $p((z_0, w_0))$. Then ${\mathcal M}_0 = \bigcup_n {\mathcal M}_0^n$ and, hence, also $\mathfrak{D} = \bigcup_n \mathfrak{D}_n$, where $\mathfrak{D}_n \coloneqq p_u({\mathcal M}_0^n)$ and $p_u : {\mathbb C}_z \times {\mathbb R}_u \to {\mathbb C}_z$ is the canonical projection. Since the projection of ${\mathcal L}_0$ onto ${\mathcal M}_0$ is one-to-one, we see that the function $U_n : \mathfrak{D}_n \to {\mathbb R}_u$ such that its graph $\Gamma(U_n)$ coincides with ${\mathcal M}_0^n$ has a single-valued conjugate function. Then, by Lemma 3, we know that the domain $\mathfrak{D}_n$ has to be simply-connected. It follows now that the domain $\mathfrak{D}$, as the union of an increasing sequence of simply-connected domains, is also simply-connected.

For proving that $\mathfrak{D} = {\mathbb C}_z$ we need to show that $\partial\mathfrak{D} = \emptyset$. Assume to the contrary that $\partial\mathfrak{D} \neq \emptyset$. In this case we conclude from simply-connectedness of the domain $\mathfrak{D}$ that $\partial\mathfrak{D}$ has the full harmonic measure relatively $\mathfrak{D}$. Then, since the function $\log|{\mathcal H}(z)|$ is harmonic near $\partial\mathfrak{D}$, it follows from the established in Lemma 11 property  $\log|{\mathcal H}(z)| \to + \infty$ as $z \to \partial{\mathfrak{D}}$ that this function has to be identically equal to $+ \infty$. This gives the desired contradiction and completes the proof of Lemma 12.
\end{proof}

The last lemma of this section gives us the required foliation structure of the boundary $\partial{\mathfrak{A}}$ of $\mathfrak{A}$.

\vspace{3mm}
\noindent \textbf{Lemma 13.} {\it Let $\mathfrak{A}$ be a pseudoconvex domain of the form 
\[ \mathfrak{A} = \{(z, w) \in \mathbb C^2 : v > F(z, u)\}, \]
where $w = u + iv$ and F is a continuous function on ${\mathbb C}_z \times {\mathbb R}_u$. Assume that there exists an entire function ${\mathcal H} : {\mathbb C}_z \to {\mathbb C}_w$ such that for its graph $\Gamma({\mathcal H})$ the inclusion $\Gamma({\mathcal H}) \subset \overline{\mathfrak{A}}$ holds. Then the boundary $\partial{\mathfrak{A}}$ of $\mathfrak{A}$ can be represented as 
\[ \partial{\mathfrak{A}} = \bigcup_{t \in {\mathbb R}} \Gamma({\mathcal H} + (t +i\Theta(t))), \]
where ${\mathcal H} : {\mathbb C}_z \to {\mathbb C}_w$  is an entire function, $\Theta : {\mathbb R} \to {\mathbb R}$ is a continuous function and $\Gamma({\mathcal H} + (t +i\Theta(t)))$ denote the graph of the function $w = {\mathcal H}(z) + (t +i\Theta(t))$.

In particular, every continuous Levi flat graph over ${\mathbb C}_z \times {\mathbb R}_u$ is foliated by translations in the $w$-direction of the graph of just one entire function.} 

\begin{proof} Observe first that the inequality defining the Model domain $\mathfrak{A}$ implies that for the function $w = {\mathcal H}(z) + i$ we have the inclusion $\Gamma({\mathcal H} + i) \subset \mathfrak{A}$. Consider a biholomorphic change of coordinates ${\mathfrak F}_1$ in $\mathbb{C}_{z, w}^2$:

\[ z \to z \eqqcolon z', \, w \to w - ({\mathcal H}(z) + i) \eqqcolon w'. \]

\noindent
Note that ${\mathfrak{A}}_1 \coloneqq {\mathfrak F}_1(\mathfrak{A})$ is also a Model domain and, hence, the inclusion 

\[ {\mathbb{C}_z} \times \{it, t \in [0, + \infty) \} \subset {\mathfrak F}_1(\mathfrak{A}) \]
holds. Now we perform one more biholomorphic transformation 

\[ {\mathfrak F}_2 : {\mathbb{C}_z} \times ({\mathbb{C}_w} \setminus \{it, t \in [0, + \infty) \}) \to {\mathbb{C}_z} \times \Delta_1(0), \]

\noindent
namely,

\[ z' \to z' \eqqcolon z'', \, w' \to \frac{\sqrt{-iw'} - i}{\sqrt{-iw'} + i} \eqqcolon w'', \]
where in the square root we choose the branch that maps the set $\mathbb{C} \setminus \{it, t \in [0, + \infty) \}$ to the upper half-plane. Observe that the set ${\mathfrak F}_2({\mathfrak F}_1(\mathbb{C}_{z, w}^2 \setminus \mathfrak{A})) \cup ({\mathbb{C}_z} \times \{1\})$ will be closed and pseudoconcave. Moreover, it is contained in ${\mathbb{C}_z} \times \overline{\Delta}_1(0)$, i.e. this set is an analytic multifunction which is also  bounded. Then Liouville's theorem for analytic multifunctions (see, for example, Theorem 1.11 in \cite{Aupetit94}) tells us that this set has the form ${\mathbb{C}_z} \times E$ for some compact set $E \subset \overline{\Delta}_1(0)$. Therefore, since the boundary of the set ${\mathfrak F}_2({\mathfrak F}_1(\mathbb{C}_{z, w}^2 \setminus \mathfrak{A})) \cup ({\mathbb{C}_z} \times \{1\}) = {\mathbb{C}_z} \times E$ is the disjoint union of lines parallel to the z-axis, we conclude, in view of the form of the maps ${\mathfrak F}_1$ and ${\mathfrak F}_2$, that the boundary of the set $\mathbb{C}_{z, w}^2 \setminus \mathfrak{A} = {\mathfrak F}_1^{-1}(({\mathfrak F}_2^{-1}({\mathbb{C}_z} \times E))$ is constituted by the translations in the w-direction of the graph $\Gamma({\mathcal H})$ of the function ${\mathcal H}$. This completes the proof of Lemma 13.
\end{proof}

Now we can finally conclude from the construction of the limit leaf ${\mathcal L}_0$ and the properties established in Lemmas 11 -$\,$13 that the following main statement of this section holds true.

\vspace{3mm}
\noindent \textbf{Proposition 2.} {\it If $\, \widetilde{\mathfrak{A}} \neq \mathfrak{A}$, then 
there  is an entire function ${\mathcal H} : {\mathbb C}_z \to {\mathbb C}_w$ and a continuous function $\Theta : {\mathbb R} \to {\mathbb R}$ such that the boundary $\partial{\mathfrak{A}}$ of $\mathfrak{A}$ can be represented as
\[ \partial{\mathfrak{A}} = \bigcup_{t \in {\mathbb R}} \Gamma({\mathcal H} + (t +i\Theta(t))), \]
where $\Gamma({\mathcal H} + (t +i\Theta(t)))$ denote the graph of the function $w = {\mathcal H}(z) + (t +i\Theta(t))$.}

\section{Proof of the Main Theorem}

\noindent
Now we can finally complete the proof of the Main Theorem.

\begin{proof}
Observe first that, since, in the case when $\widetilde{\mathfrak{A}} = \mathfrak{A}$, we know by Proposition 1 that $\mathfrak{A}$ possesses a bounded smooth strictly plurisubharmonic function, it follows from the definition of the core that $\mathfrak{c}(\mathfrak{A}) = \emptyset$. It also follows from Sibony's Theorem 2 above that in this case the domain $\mathfrak{A}$ is Kobayashi hyperbolic, and from Chen-Zhang's Theorem 3 above that $\mathfrak{A}$ possesses a Bergman metric.

Since, in the second case when $\widetilde{\mathfrak{A}} \neq \mathfrak{A}$, we have by Proposition 2 the required foliation of the boundary $\partial{\mathfrak{A}}$  of the domain $\mathfrak{A}$ by translations of the graph of an entire function, we conclude that to finish the proof of the theorem we only need to show that in this case:
\begin{enumerate}
  \item[(i)] $\mathfrak{c}(\mathfrak{A}) \neq \emptyset$.

  \item[(ii)] The domain $\mathfrak{A}$ is not Brody hyperbolic.

  \item[(iii)] The domain $\mathfrak{A}$ is not Kobayashi hyperbolic.

  \item[(iv)] The domain $\mathfrak{A}$ does not possess a Bergman metric.
\end{enumerate}

For showing (i), (ii) and (iii) it is enough to observe that the graph of the entire function $w = {\mathcal H}(z) + i$ is contained in $\mathfrak{A}$. Then, using Liouville's theorem on this curve, we see that $\Gamma({\mathcal H} + i) \subset \mathfrak{c}(\mathfrak{A})$ and hence $\mathfrak{c}(\mathfrak{A}) \neq \emptyset$. We also see that the domain $\mathfrak{A}$, by the definition, is not Brody hyperbolic, and that the Kobayashi distance along the curve $\Gamma({\mathcal H} + i)$ degenerates and, hence, in this case the domain $\mathfrak{A}$ can not be Kobayashi hyperbolic.

For proving (iv) we can proceed the same way as in the proof of Proposition 2. Namely, we first perform a biholomorphic change of coordinates ${\mathfrak F}_1$ in $\mathbb{C}_{z, w}^2$:

\[ z \to z \eqqcolon z', \, w \to w - {\mathcal H}(z) \eqqcolon w'. \]
Then the domain $\mathfrak {F}_1(\mathfrak{A})$ will have the form ${\mathfrak F}_1(\mathfrak{A}) = \{(z', w') \in \mathbb C^2 : v' > \Theta(u')\}$. Since the domain $\Sigma \coloneqq \{ w' \in \mathbb C : v' > \Theta(u') \}$, where $w' = u' + iv'$, is simply-connected, there exists a conformal map  ${\mathfrak f} : \Sigma \to \Delta_1(0)$ to the unit disk $\Delta_1(0)$ in $\mathbb{C}_{w'}$. Hence, if we consider the map $\mathfrak {F}_2 : {\mathbb C}_z \times \Sigma \to {\mathbb C}_z \times \Delta_1(0)$ defined as:

\[ z' \to z' \eqqcolon z'', \, w' \to {\mathfrak f}(w') \eqqcolon w'', \]
then $\mathfrak {F}_2  \circ \mathfrak {F}_1$ will be a biholomorphic map of the domain $\mathfrak{A}$ onto the unbounded cylinder $\mathbb{C}_z \times \Delta_1(0)$. But the space of holomorphic $L^2$-functions on this cylinder consists of just one function which is identically equal to zero. Therefore the space of holomorphic $L^2$-functions on $\mathfrak{A}$ also consists of just one function which is equal to zero. It follows then that the domain $\mathfrak{A}$ can not have a Bergman metric. The proof of the Main Theorem is now completed.
\end{proof}

\vspace{3mm}
\noindent \textbf{Remark 2.} In the case when the boundary $\partial{\mathfrak A}$ of the domain $\mathfrak A$ is not foliated by translations of the graph of an entire function, we know from our Main Theorem that its core ${\mathfrak c}(\mathfrak A)$ is empty. It follows then from Theorem 1 in \cite{GallagherHarzHerbort17} that in this case $\dim{A^2(\mathfrak A)} = \infty$, where $A^2(\mathfrak A)$ is the Bergman space, i.e. the space of holomorphic functions on $\mathfrak A$ which belong to $L^2(\mathfrak A)$.

\vspace{3mm}
\noindent \textbf{Remark 3.} Since, by Lemma 13, every continuous Levi flat graph over ${\mathbb C}_z \times {\mathbb R}_u$ is foliated by translations of the graph of just one entire function, and since, by Theorem 1 above, the Dirichlet problem for Levi flat graphs over bounded convex domains is solvable for any continuous boundary function, we conclude that not all Levi flat graphs over bounded convex domains can be uniformly approximated by the globally defined (i.e. defined over the whole of ${\mathbb C}_z \times {\mathbb R}_u$) Levi flat graphs.

\vspace{3mm}
\noindent \textbf{Remark 4.} We do not know if all the Model domains $\mathfrak A \subset {\mathbb C}^2$, whose boundary $\partial{\mathfrak A}$ is not foliated by translations of the graph of an entire function, possess a Carath$\acute{e}$odory metric.

\vspace{3mm}
\noindent
{\bf Acknowledgements.} {\it Part of this work was done while the author was a visitor at the Capital Normal University (Beijing). It is his pleasure to thank this institution for its hospitality and good working conditions. The author also would like to thank B.-Y. Chen, K.-T. Kim and W. Zwonek for useful information related to Remark 2, K.-T. Kim for references on some papers where Model domains were previously studied and A. Isaev for references related to the Lie group structure on ${\text{Aut}}(\mathfrak{A})$ for Kobayashi hyperbolic domains $\mathfrak{A}$.}

%
%
 \vspace{1truecm}

%
%
%
{\sc N. Shcherbina: Department of Mathematics, University of Wuppertal --- 42119 Wuppertal, Germany}
  
{\em e-mail address}: {\texttt shcherbina@math.uni-wuppertal.de}

\end{document}